\documentclass[10pt,twocolumn]{article}

\usepackage{mathrsfs}
\usepackage{abstract}

\usepackage{enumitem}
\setenumerate[0]{leftmargin=*}

\usepackage{graphicx}
\usepackage{multirow,wrapfig}
\usepackage{amsmath, amsthm, amsfonts, amssymb}
\usepackage[linesnumbered, ruled, vlined]{algorithm2e}
\usepackage{xcolor,natbib}
\usepackage[normalem]{ulem}
\usepackage{makecell}

\usepackage{amssymb} 

\usepackage[unicode=true, bookmarks=false, 
 breaklinks=false,pdfborder={0 0 1}]
 {hyperref}
\hypersetup{
 colorlinks,citecolor=blue,filecolor=blue,linkcolor=blue,urlcolor=blue}
 \usepackage{url}

\newtheorem{theorem}{Theorem}
\newtheorem{lemma}{Lemma}

\allowdisplaybreaks

\definecolor{yc}{RGB}{190,0,255}  
\definecolor{dacong}{RGB}{10,103,68}

\newcommand{\ex}[2]{\mathbb{E}_{#1}\left[#2\right]}

\newcommand{\prnbig}[1]{\big({#1}\big)} % parentheses
\newcommand{\norm}[1]{\left\|{#1}\right\|} % norm
\newcommand{\innprod}[1]{\big\langle{#1} \big\rangle} % inner product
\newcommand{\cS}{\mathcal{S}}
\newcommand{\cA}{\mathcal{A}}
\newcommand{\cB}{\mathcal{B}}
\newcommand{\cO}{\mathcal{O}}

\newcommand{\cH}{\mathcal{H}}

\newcommand{\KL}[2]{\mathsf{KL}\prnbig{{#1}\,\|\,{#2}}}

\title{Global Convergence of Policy Gradient Methods \\ in Reinforcement Learning, Games and Control}
 	
\author{Shicong Cen \\
Carnegie Mellon University
\and Yuejie Chi\\
Carnegie Mellon University\thanks{S. Cen and Y. Chi are with the Department of Electrical and Computer Engineering, Carnegie Mellon University; emails: {\tt\small \{shicongc,yuejiec\}@andrew.cmu.edu}.}  
}

\date{}

\begin{document}

\twocolumn[ 
\maketitle
 
\begin{abstract}
Policy gradient methods, where one searches for the policy of interest by maximizing the value functions using first-order information, become increasingly popular for sequential decision making in reinforcement learning, games, and control. Guaranteeing the global optimality of policy gradient methods, however, is highly nontrivial due to nonconcavity of the value functions. In this exposition, we highlight recent progresses in understanding and developing policy gradient methods with global convergence guarantees, putting an emphasis on their finite-time convergence rates with regard to salient problem parameters.

\end{abstract}

\bigskip
]

\saythanks

\section{Introduction} 

Sequential decision making is a canonical task that lies at the heart of a wide spectrum of disciplines such as reinforcement learning (RL), games and control, finding numerous applications in autonomous driving, robotics, supply chain management, resource scheduling, and so on. While classical approaches such as dynamic programming \citep{bertsekas2017dynamic} require model knowledge, modern practices advocate for the model-free approach as an alternative: one seeks the policy of interest directly based on data collected through interactions with the environments, without estimating the model. The advantage is that the model-free approach is often more memory efficient, as well as more agile to changes.

One  prevalent class of model-free approaches is policy gradient (PG) methods, which follow an optimizer's perspective by formulating a value maximization problem with regard to parameterized policies, and performing gradient updates to improve the policy iteratively --- often based on noisy feedbacks received from the environments. PG methods and their variants have become the de facto standard practices in an increasing number of domains, due to their seamless integration with neural network parameterization and adaptivity to various problem setups involving discrete, continuous or mixed action and state spaces.

Despite the great empirical success of PG methods, little is known about their theoretical convergence properties --- especially when it comes to finite-time global convergence --- due to the notorious nonconcavity of the value functions. Until very recently, such understandings begin to emerge, leveraging the fact that the PG methods are typically operated on highly structured model classes, whose induced optimization landscapes turn out to be much more benign and tractable than previously thought. The execution of  problem-dependent tailored analyses, rather than relying on black-box optimization theory, fuels recent breakthroughs pioneered by \citet{fazel2018global,agarwal2019optimality,bhandari2019global}, to name just a few.
The purpose of this article is to survey the latest efforts in understanding the global convergence of PG methods in the fields of RL, game theory and control, as well as highlight algorithmic ideas that enable fast global convergence, especially with regard to salient problem parameters that are crucial in practice.

\paragraph{Organization.}
The rest of this article is organized as follows. Section~\ref{sec:RL} reviews PG methods in single-agent RL, focusing on solving tabular Markov decision processes (MDP). Section~\ref{sec:games} reviews PG methods in the game and multi-agent RL setting, using the two-player zero-sum matrix game and two-player zero-sum Markov game as illustrative examples. Section~\ref{sec:control} moves onto PG methods in control, focusing on solving the linear quadratic regulator (LQR). We conclude in Section~\ref{sec:conclusions}. 

\paragraph{Notation.} We use $\|A\|$, $\|A\|_{\mathsf{F}}$ and $\sigma_{\min}(A)$ to represent the spectral norm, the Frobenius norm, and the smallest singular value of a  matrix $A$, respectively. For two vectors $a$ and $b$, we use $\tfrac{a}{b}$ to denote their entrywise division when exists, and $\|a\|_\infty$ denotes the entrywise maximum absolute value of a vector $a$. Given a set $\cS$, we let $\Delta(\cS)$ represent the probability simplex over $\cS$. Given two distributions $p$ and $q$, $\KL{p}{q}$ denotes the Kullback-Leibler (KL) divergence from $q$ to $p$. Last but not least, let $\mathcal{P}_{\mathcal{C}}(\cdot)$ be the projection operator onto the set $\mathcal{C}$. To characterize the iteration complexity, we write $T(x) = \cO(f(x))$ when $T(x) \le C f(x), \forall x \ge X$ for some $C, X > 0$. Here $x$ is typically set to $1/\epsilon$ or other salient problem parameters. We use $\widetilde{\cO}(\cdot)$ to suppress logarithmic factors from the standard order notation $\cO(\cdot)$.

\section{Global convergence of policy gradient methods in RL}
\label{sec:RL}

In this section, we review recent progresses in developing policy gradient methods for single-agent RL, focusing on the basic model of tabular Markov decision processes (MDPs).

\subsection{Problem settings}

\paragraph{Markov decision processes.} An infinite-horizon discounted MDP models the sequential decision making problem as $\mathcal{M} = (\cS, \mathcal{A}, P, r, \gamma)$ with the state space $\cS$, the action space $\mathcal{A}$, the transition kernel $P: \cS\times\mathcal{A} \to \Delta(\cS)$, the reward function $r:\cS\times\mathcal{A}\to [0,1]$ and the discount factor $\gamma \in (0,1)$. Upon the agent choosing action $a\in \mathcal{A}$ at state $s \in \cS$, the environment will move to a new state $s'$ according to the transition probability $P(s'|s, a)$ and assign a reward $r(s,a)$ to the agent. The action selection rule is implemented by a randomized policy $\pi:\cS\to\Delta(\mathcal{A})$, where $\pi(a|s)$ specifies the probability of choosing action $a$ in state $s$. {We shall focus exclusively on the case where both $\mathcal{S}$ and $\mathcal{A}$ are finite throughout this article.}

\paragraph{Value functions and Q-functions.} The value function $V^\pi: \cS \to \mathbb{R}$ is defined as the expected discounted cumulative reward starting at state $s$:
\begin{equation*}
    V^\pi(s) := \mathbb{E}\bigg[\sum_{t=0}^\infty \gamma^t r(s_t, a_t)\,\big|\,s_0 = s\bigg],
\end{equation*}
where $a_t \sim \pi(\cdot|s_t)$ is obtained by executing policy $\pi$ and $s_{t+1} \sim P(\cdot|s_t, a_t)$ is generated by the MDP. The Q-function (or action-value function) $Q^\pi: \cS\times\mathcal{A}\to \mathbb{R}$ is defined in a similar manner with an initial state-action pair $(s, a)$:
\begin{equation*}
    Q^\pi(s,a) := \mathbb{E}\bigg[\sum_{t=0}^\infty \gamma^t r(s_t, a_t)\,\big|\,s_0 = s, a_0 = a\bigg].
\end{equation*}
In addition, the advantage function of policy $\pi$ is defined as
\begin{equation*}
    A^\pi(s, a) := Q^\pi(s,a) - V^\pi(s).
\end{equation*}
It is well-known that there exists an optimal policy $\pi^\star$ that maximizes the value function $V^\pi(s)$ for all $s\in\cS$ simultaneously \citep{bellman1952theory,puterman2014markov}. We denote the resulting optimal value function and $Q$-function by $V^\star$ and $Q^\star$, which satisfy the well-known Bellman optimality equations. The optimal policy $\pi^{\star}$ can be implied from the optimal $Q$-function in a greedy fashion as
\begin{equation} \label{eq:greedy}
 \pi^{\star}(a|s) = \arg\max_{a\in\cA}\; Q^{\star}(s,a), \qquad \forall \;s \in\cS.
\end{equation}

\paragraph{Policy parameterizations.} Given some prescribed initial distribution $\rho$ over $\cS$, policy optimization methods seek to maximize the value function
\begin{equation*}
 V^{\pi}(\rho) := V^{\pi_\theta}(\rho)
\end{equation*}
over the policy $\pi$, where $\pi: = \pi_\theta$ is often parameterized via some parameter  $\theta$, and we overload the notation to denote by $V^\pi(\rho)$ the expected value function $\mathbb{E}_{s\sim \rho}\big[V^\pi(s)\big]$. Note that it is straightforward to observe that $V^\star(\phi) - V^{\pi}(\phi) \le \|\phi/\rho\|_\infty(V^\star(\rho) - V^{\pi}(\rho))$ for general choices of $\phi \in \Delta(\cS)$, which characterizes the effect of possible discrepancy between the initial distributions in training and deploying. Common choices of policy parameterization are as follows.
\begin{itemize}
    \item \textit{Direct parameterization:} the policy is directly parameterized by
    \begin{equation*}
        \pi_\theta(a|s) = \theta(s,a),
    \end{equation*}
    where  $\theta \in  \{\theta \in \mathbb{R}^{|\cS||\mathcal{A}|}:  \, \, \theta(s,a)\ge 0, \sum_{a\in\mathcal{A}}\theta(s,a)=1\}$.
    \item \textit{Tabular softmax parameterization:} For $\theta \in \mathbb{R}^{|\cS||\mathcal{A}|}$, the policy $\pi_\theta$ is generated through softmax transform
    \begin{equation*}
        \pi_\theta(a|s) = \frac{\exp(\theta(s,a))}{\sum_{a'\in\mathcal{A}}\exp(\theta(s,a'))},
    \end{equation*}
    leading to an unconstrained optimization over $\theta$.  Although beyond the scope of this exposition, softmax prameterization is more popular with function approximation, by replacing $\theta(s,a)$ with $f_\theta(s,a)$ typically implemented by neural networks. 
\end{itemize}

\paragraph{Policy gradients.} The gradient $\nabla_\theta V^{\pi_\theta}(\rho)$ plays an instrumental role in developing first-order policy optimization methods. To facilitate presentation, we shall first introduce the discounted state visitation distributions of a policy $\pi$:
\begin{equation*}
    d_{s_0}^\pi(s) := (1-\gamma) \sum_{t=0}^\infty \gamma^t \mathbb{P}(s_t = s|s_0),
\end{equation*}
where the expectation is with respect to the trajectory $(s_0, a_0, s_1, a_1, \cdots)$ generated by the MDP under policy $\pi$. We further denote by $d_\rho^\pi$ the discounted state visitation distribution when $s_0$ is randomly drawn from distribution $\rho$, i.e.,
\begin{equation*}
    d_\rho^\pi(s) := \mathbb{E}_{s_0\sim \rho}\big[d_{s_0}^\pi(s)\big].
\end{equation*} 

The policy gradient of parameterized policy $\pi_\theta$ is then given by \citep{williams1992simple}
\begin{align*}
    \nabla_\theta V^{\pi_\theta}(\rho) & = \frac{1}{1-\gamma}\mathbb{E}_{s\sim d_\rho^{\pi_\theta}, a\sim \pi_\theta(\cdot|s)}\big[\nabla_\theta \log \pi_\theta (a|s) Q^{\pi_\theta}(s,a)\big]  \\
     & = \frac{1}{1-\gamma}\mathbb{E}_{s\sim d_\rho^{\pi_\theta}, a\sim \pi_\theta(\cdot|s)}\big[\nabla_\theta \log \pi_\theta (a|s) A^{\pi_\theta}(s,a)\big]  ,
\end{align*}
which can be evaluated, for example, via REINFORCE \citep{williams1992simple}. The use of the advantage function $A^{\pi_\theta}(s,a)$, rather than the Q-function $Q^{\pi_\theta}(s,a)$, often helps to reduce the variance of the estimated policy gradient.

For notational simplicity, we shall denote by $\theta^{(t)}$ and $\pi^{(t)}$ the parameter and the policy at the $t$-th iteration, and use $V^{(t)}$, $Q^{(t)}$, $A^{(t)}$, $d_\rho^{(t)}$ to denote $V^{\pi^{(t)}}$, $Q^{\pi^{(t)}}$, $A^{\pi^{(t)}}$, $d_\rho^{\pi^{(t)}}$, respectively. In addition, we assume the policy gradients and the value functions are exactly evaluated throughout this article, which enables us to focus on the optimization aspect of PG methods.  

\subsection{Projected policy gradient method}
The most straightforward first-order policy optimization method is to adopt direct parameterization and perform projected gradient ascent updates:
\begin{equation}
    \theta^{(t+1)} = \mathcal{P}_{\Delta(\mathcal{A})^{|\cS|}}\big(\theta^{(t)} + \eta \nabla_\theta V^{(t)}(\rho) \big),
    \label{eq:direct_pg}
\end{equation}
or equivalently
\begin{equation*}
    \pi^{(t+1)} = \mathcal{P}_{\Delta(\mathcal{A})^{|\cS|}} \big(\pi^{(t)} + \eta \nabla_\pi V^{(t)}(\rho) \big),
    % \label{eq:direct_pg}
\end{equation*}
where $\eta>0$ is the learning rate, and
\begin{equation*}
    \nabla_{\theta(s,a)} V^{(t)}(\rho) = \nabla_{\pi(s,a)} V^{(t)}(\rho) = \frac{1}{1-\gamma} d_\rho^{(t)}(s)Q^{(t)}(s,a).
\end{equation*}
As the value function $V^{\pi_\theta}(\rho)$ is $\frac{2\gamma|\mathcal{A}|}{(1-\gamma)^3}$-smooth \citep{agarwal2019optimality}, setting the learning rate to $0 < \eta \le \frac{(1-\gamma)^3}{2\gamma|\mathcal{A}|}$ ensures monotonicity of $V^{(t)}(\rho)$ in  $t$. On the other end, it is critically established in \citet{agarwal2019optimality} that the value function satisfies the following gradient domination condition.
\begin{lemma}[Variational gradient domination]
    For any policy $\pi$, we have
    \begin{equation*}
        V^\star(\rho) - V^{\pi}(\rho) \le \frac{1}{1-\gamma}\bigg\|\frac{d_\rho^{\pi^{\star}}}{\rho}\bigg\|_\infty \max_{\pi'\in\Delta(\mathcal{A})^{|\cS|}}(\pi' - \pi)^\top \nabla_\pi V^{\pi}(\rho).
    \end{equation*}    
\end{lemma}
The above lemma associates the optimality gap $V^\star(\rho) - V^{\pi}(\rho)$ with a variational gradient term, allowing the iterates to converge globally as stated below.

\begin{theorem}[\mbox{\cite{agarwal2019optimality}}] With $0 < \eta \le \frac{(1-\gamma)^3}{2\gamma|\mathcal{A}|}$, the iterates of the projected PG method \eqref{eq:direct_pg} satisfies
\begin{equation*}
    \min_{0\le t \le T} V^\star (\rho) - V^{(t)}(\rho) \le \frac{4\sqrt{|\cS|}}{1-\gamma}\bigg\|\frac{d_\rho^{\pi^{\star}}}{\rho}\bigg\|_\infty \sqrt{\frac{2(V^\star (\rho) - V^{(0)}(\rho))}{\eta T}}.
\end{equation*}
\label{thm:agarwal_pgd}
\end{theorem}
Theorem \ref{thm:agarwal_pgd} establishes an iteration complexity of $\mathcal{O}\Big(\frac{|\cS||\mathcal{A}|}{(1-\gamma)^6\epsilon^2}\big\|\frac{d_\rho^{\pi^{\star}}}{\rho}\big\|_\infty^2\Big)$ for finding an $\epsilon$-optimal policy, which is later improved to $\mathcal{O}\Big(\frac{|\cS||\mathcal{A}|}{(1-\gamma)^5\epsilon}\big\|\frac{d_\rho^{\pi^{\star}}}{\rho}\big\|_\infty^2\Big)$ \citep{xiao2022convergence}. However, the projection operator introduces $\mathcal{O}(\log |\mathcal{A}|)$ computational overhead every iteration and is less compatible with  function approximation. This motivates the study of PG methods that are compatible with unconstrained optimization, e.g., by using softmax parameterization.

\subsection{Softmax policy gradient method}
With softmax parameterization, the policy gradient method writes
\begin{equation}
    \theta^{(t+1)} = \theta^{(t)} + \eta\nabla_\theta V^{(t)}(\rho),
    \label{eq:softmax_pg}
\end{equation}
where
\begin{equation*}
    \nabla_{\theta(s,a)} V^{(t)}(\rho) = \frac{\eta}{1-\gamma}d_\rho^{(t)}(s)\pi^{(t)}(a|s)A^{(t)}(s,a).
\end{equation*}
Remarkably, \citet{agarwal2019optimality} established the asymptotic global convergence of the softmax PG method as follows.
\begin{theorem}[\mbox{\cite{agarwal2019optimality}}]
    With constant learning rate $0 < \eta \le (1-\gamma)^3/8$, the softmax PG method converges to the optimal policy, i.e., $V^{(t)}(s) \to V^\star(s)$ as $t\to \infty$ for all $s \in \cS$.
\end{theorem}

\citet{mei2020global} later demonstrated an iteration complexity of $\mathcal{O}(\frac{1}{c(\mathcal{M})^2\epsilon})$ for achieving an $\epsilon$-optimal policy, where $c(\mathcal{M})$ is a trajectory-dependent quantity depending on salient problem parameters such as the number of states $|\cS|$ and the effective horizon $(1-\gamma)^{-1}$. Unfortunately, this quantity $c(\mathcal{M})$ can be rather small and does not exclude the possibility of incurring excessively large iteration complexity, as demonstrated by the following hardness result \citep{li2021softmax}.
\begin{theorem}[\cite{li2021softmax}]
    There exist universal constants $c_1, c_2, c_3 > 0$ such that for any $\gamma \in (0.96, 1)$ and $|\cS| \ge c_3 (1-\gamma)^{-6}$, one can find a $\gamma$-discounted MDP such that the softmax PG method takes at least
    \begin{equation*}
        \frac{c_1}{\eta}|\cS|^{2^{\frac{c_2}{1-\gamma}}}
    \end{equation*}
    iterations to reach $\|V^\star - V^{(t)}\|_\infty \le 0.15$.
\end{theorem}
Therefore, though guaranteed to converge globally, softmax PG can take (super-)exponential time to even reduce the optimality gap within a constant level. 
Intuitively speaking, softmax PG method fails to achieve a reasonable convergence rate when the probability $\pi^{(t)}(a^\star(s)|s)$ assigned to the optimal action $a^\star(s)$ is close to zero. \citet{agarwal2019optimality} proposed to penalize the policy for getting too close to the border, by imposing a log barrier regularization 
\begin{equation*}
    V_\omega^{\pi_\theta}(\rho) = V^{\pi_\theta}(\rho) + \frac{\omega}{|\cS||\mathcal{A}|}\sum_{s\in\cS, a\in\mathcal{A}}\log\pi_\theta(a|s).
\end{equation*}
With an appropriate choice of regularization parameter $\omega$, the regularized softmax PG can achieve an $\epsilon$-optimal policy within $\mathcal{O}\Big(\frac{|\cS|^2|\mathcal{A}|^2}{(1-\gamma)^6\epsilon^2}\Big\|\frac{d_\rho^\star}{\rho}\Big\|_\infty^2\Big)$ iterations \citep{agarwal2019optimality}. Nonetheless, this regularization scheme is not as popular in practice, compared to the entropy regularization scheme that will be discussed momentarily.

\subsection{Natural policy gradient method}
Both projected PG and softmax PG fall short of attaining an iteration complexity that is independent of salient problem parameters, especially with respect to the size of the state space $|\mathcal{S}|$. This ambitious goal can be achieved, somewhat surprisingly, by adopting Fisher information matrix as a preconditioner, which leads to the natural policy gradient (NPG) method \citep{kakade2002natural}:
\begin{equation}    \label{eq:softmax_npg}
        \theta^{(t+1)} = \theta^{(t)} + \eta (\mathcal{F}_\rho^{\theta^{(t)}})^\dagger \nabla_\theta V^{(t)}(\rho),
\end{equation}
where
$$   \mathcal{F}_\rho^\theta := \mathbb{E}_{s\sim d_\rho^{\pi_\theta}, a\sim \pi_\theta(\cdot|s)}\Big[\big(\nabla_\theta\log \pi_\theta(a|s)\big)\big(\nabla_\theta\log \pi_\theta(a|s)\big)^\top\Big]$$
is the Fisher information matrix,  and $^\dagger$ denotes the Moore-Penrose pseudoinverse. With softmax parameterization, the NPG updates take the form 
\begin{equation*}
    \theta^{(t+1)} = \theta^{(t)} + \frac{\eta}{1-\gamma}A^{(t)},
\end{equation*}
or equivalently,
\begin{equation*}
    \pi^{(t+1)}(a|s) \propto \pi^{(t)}(a|s) \exp\Big(\frac{\eta Q^{(t)}(s,a)}{1-\gamma}\Big),
\end{equation*}
It is noted that the (softmax) NPG update rule coincides with the multiplicative weights update (MWU) method \citep{cesa2006prediction}, and that the update rule does not depend on the initial state distribution $\rho$. \citet{shani2019adaptive} first established a global convergence rate of $\mathcal{O}\big(\frac{1}{(1-\gamma)^2\sqrt{T}}\big)$ using decaying learning rate $\eta_t = \mathcal{O}\big(\frac{1-\gamma}{\sqrt{t}}\big)$, which was improved by \citet{agarwal2019optimality} using a constant learning rate $\eta$, stated as follows.
\begin{theorem}[\mbox{\cite{agarwal2019optimality}}]     \label{thm:npg}
    With uniform initialization $\theta^{(0)} = 0$ and constant learning rate $\eta > 0$, the iterates of NPG satisfy
    \begin{equation*}
        V^\star(\rho) - V^{(T)}(\rho) \le \frac{1}{T}\Big(\frac{\log|\mathcal{A}|}{\eta} + \frac{1}{(1-\gamma)^2}\Big).
    \end{equation*}
     %   for $T \ge 1$.
\end{theorem}
Encouragingly, as long as $\eta\geq \frac{(1-\gamma)^2}{\log |\mathcal{A}|}$, the iteration complexity of NPG methods becomes $\mathcal{O}\big(\frac{1}{(1-\gamma)^2 T}\big)$, which is independent of the size of the state-action space. 
On the complementary side, the iteration complexity of NPG is lower bounded by $\frac{\Delta}{(1-\gamma)|\mathcal{A}|}\exp(-\eta\Delta T)$ --- established in \citet{khodadadian2021linear} --- where the optimal advantage function gap $\Delta  = \min_s\min_{a\neq a^{\star}(s)} |A^{\star}(s,a)| \ge 0$ is determined by the MDP instance. As the lower bound attains its maximum $\frac{1}{(1-\gamma)e\eta T}$ when $\Delta = \frac{1}{\eta T}$, it is immediate that the sublinear rate in Theorem \ref{thm:npg} cannot be improved in $T$. Nonetheless, two strategies to achieve even faster linear convergence with NPG updates include (i) adopting increasing/adaptive learning rates \citep{khodadadian2021linear,bhandari2020note,lan2021policy,xiao2022convergence}, or (ii) introducing entropy regularization \citep{cen2020fast,lan2021policy,zhan2021policy}, which we shall elaborate on next. 

\subsection{Entropy regularization} 
Introducing entropy regularization is a popular technique in practice to promote exploration \citep{haarnoja2017reinforcement}. Specifically,  one seeks to optimize the entropy-regularized value function defined as
\begin{equation*}
        V_\tau^{\pi}(s) = V^{\pi}(s) + \frac{\tau}{1-\gamma} \mathbb{E}_{s'\sim d_s^{\pi}}\big[ \mathcal{H} (\pi(\cdot|s'))\big],
\end{equation*}
where $\mathcal{H}(\pi(\cdot|s)) = - \sum_{a\in \mathcal{A}} \pi(a|s) \log \pi(a|s)$ is the entropy of policy $ \pi(\cdot |s) $, and $\tau > 0$ serves as the regularization parameter known as the temperature. The entropy-regularized $Q$-function is defined as
$$        Q_\tau^{\pi}(s,a) = r(s,a) + \gamma  \mathbb{E}_{s'\sim P(\cdot|s,a)} \big[ V_\tau^{\pi}(s') \big]. $$
The resulting optimal value function, $Q$-function and optimal policy are denoted by $V_\tau^\star$, $Q_\tau^\star$, and $\pi_\tau^{\star}$.
From an optimization perspective, the entropy term adds curvature to the value function and ensures that the optimal policy $\pi_\tau^\star$ is unique. Interestingly, in contrast to the greedy optimal policy for the unregularized problem in \eqref{eq:greedy}, the optimal policy of the entropy-regularized problem reflects ``bounded rationality'' in decision making, namely
$$ \pi_\tau^\star (\cdot| s)\propto \exp\left( Q_\tau^{\star}(s,\cdot) /\tau \right) .$$
It should be noted, however, adding the entropy regularization generally does not make $V_\tau^\pi(\rho)$ concave unless $\tau$ is unreasonably large.  As the entropy function is bounded by $  \log|\mathcal{A}|$, the optimal entropy-regularized policy is also guaranteed to be approximately optimal for the unregularized RL problem in the following sense:
\begin{equation*}
    V^{\pi_\tau^\star}(\rho) \ge V^\star(\rho) - \frac{\tau\log|\cA|}{1-\gamma}.
\end{equation*}
  
Motivated by its benign convergence, we consider NPG for the entropy-regularized problem:
\begin{equation*}
    \theta^{(t+1)} = \theta^{(t)} + \eta (\mathcal{F}_\rho^{\theta^{(t)}})^\dagger \nabla_\theta V_\tau^{(t)}(\rho),
\end{equation*}
which can be equivalently written as
\begin{equation}     \label{eq:ent_NPG}
    \pi^{(t+1)}(a|s) \propto \pi^{(t)}(a|s)^{1-\frac{\eta\tau}{1-\gamma}} \exp\Big(\frac{\eta Q_\tau^{(t)}(s,a)}{1-\gamma}\Big).
\end{equation}
 
The following theorem shows that with appropriate choices of constant learning rate $\eta$, entropy-regularized NPG converges to the unique optimal policy $\pi_\tau^\star$ at a linear rate.
\begin{theorem} 
    For constant learning rate $0 < \eta \le (1-\gamma)/\tau$ and uniform initialization, the entropy-regularized NPG updates \eqref{eq:ent_NPG} satisfy
    \begin{equation*}
        \|V_\tau^\star - V_\tau^{(T)}\|_\infty \le \frac{15(1 + \tau \log|\mathcal{A}|)}{1-\gamma} (1-\eta\tau)^{T-1}
    \end{equation*}
and 
\ifdefined\arxiv
    \begin{align*}
        V_\tau^\star(\rho) - V_\tau^{(T)}(\rho) & \le \bigg\|\frac{\rho}{\nu_\tau^\star}\bigg\|_\infty\bigg(\frac{1+\tau\log|\mathcal{A}|}{1-\gamma}+\frac{(1-\gamma)\log|\mathcal{A}|}{\eta}\bigg)     \max\Big\{\gamma, 1-\frac{\eta\tau}{1-\gamma}\Big\}^T.
    \end{align*}
 \else
      \begin{align*}
        V_\tau^\star(\rho) - V_\tau^{(T)}(\rho) & \le \bigg\|\frac{\rho}{\nu_\tau^\star}\bigg\|_\infty\bigg(\frac{1+\tau\log|\mathcal{A}|}{1-\gamma}+\frac{(1-\gamma)\log|\mathcal{A}|}{\eta}\bigg) \\
       & \qquad \qquad\qquad \cdot  \max\Big\{\gamma, 1-\frac{\eta\tau}{1-\gamma}\Big\}^T.
    \end{align*}
    \fi
 Here, $\nu_\tau^\star$ is the stationary state distribution of policy $\pi_\tau^\star$.
\end{theorem}
The first and the second bounds are due to \citet{cen2020fast} and \citet{lan2021policy}\footnote{We discard some of the simplification steps therein and state the convergence result for a wider range of learning rate $\eta$.} respectively, where they lead to slightly different iteration complexities. Taken collectively, entropy-regularized NPG takes no more than
\begin{equation*}
    \widetilde{\cO}\Big(\min\Big\{\frac{1}{\eta\tau}\log\frac{1}{\epsilon}, \max\Big\{\frac{1}{1-\gamma}, \frac{1-\gamma}{\eta\tau}\Big\}\log\frac{\norm{{\rho}/{\nu_\tau^\star}}_\infty}{\epsilon}\Big\}\Big) 
\end{equation*}
iterations to find a policy satisfying $V_\tau^\star(\rho) - V_\tau^\pi(\rho) \le \epsilon$.
We make note that the difference stems from different analysis approaches: \citet{cen2020fast} built their analysis upon the contraction property of the soft Bellman operator (the entropy-regularized counterpart of the original Bellman operator), while \citet{lan2021policy} made use of the connection between regularized NPG and regularized mirror descent. This can be observed from the following equivalence: the update rule \eqref{eq:ent_NPG} can be equivalently expressed as
\ifdefined\arxiv
\begin{align*}
     \pi^{(t+1)}(\cdot|s)  &= \arg\min_{p\in\Delta(\cA)} \innprod{p, -Q_\tau^{(t)}(s,\cdot)} - \tau \cH(p)   + \frac{1}{\eta_{\texttt{MD}}}\KL{p}{\pi^{(t)}(\cdot|s)},
\end{align*}
\else
\begin{align*}
     \pi^{(t+1)}(\cdot|s)  &= \arg\min_{p\in\Delta(\cA)} \innprod{p, -Q_\tau^{(t)}(s,\cdot)} - \tau \cH(p) \\
    &   \qquad\qquad\qquad\qquad  + \frac{1}{\eta_{\texttt{MD}}}\KL{p}{\pi^{(t)}(\cdot|s)},
\end{align*}
\fi
with $\eta_{\texttt{MD}} = \frac{\eta}{1-\gamma-\eta\tau}$. The analysis of regularized RL can be further generalized to adopt non-strongly convex regularizers \citep{lan2021policy}, non-smooth regularizers \citep{zhan2021policy},  state-wise policy updates \citep{lan2023block}, and so on. 

\begin{figure}[t]
\begin{center}
\begin{tabular}{cc}
    \includegraphics[width=0.44\linewidth]{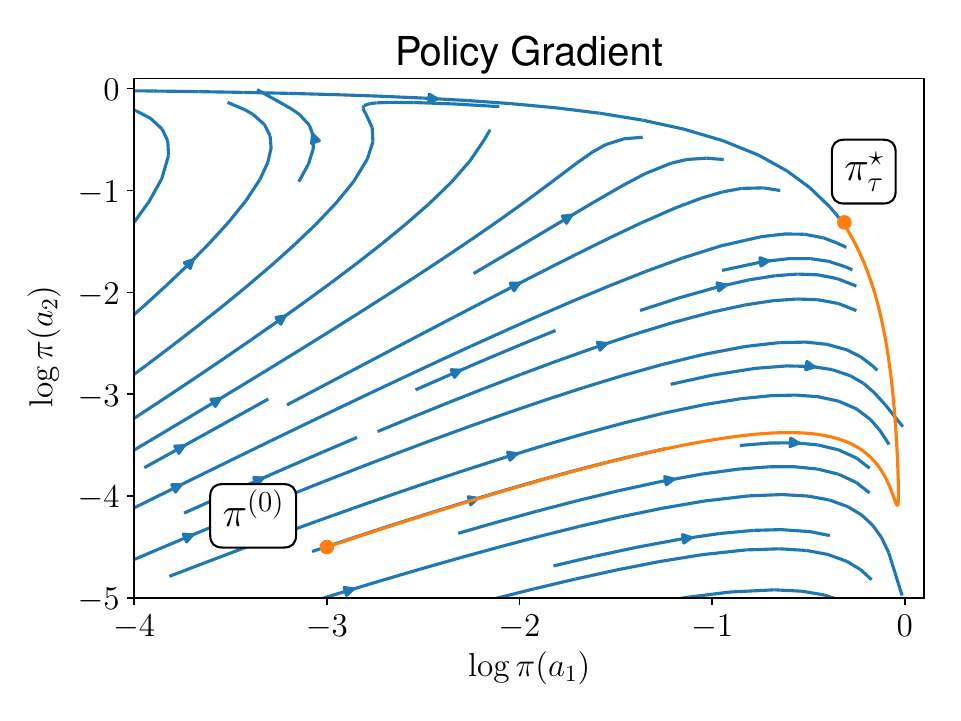} &
    \includegraphics[width=0.44\linewidth]{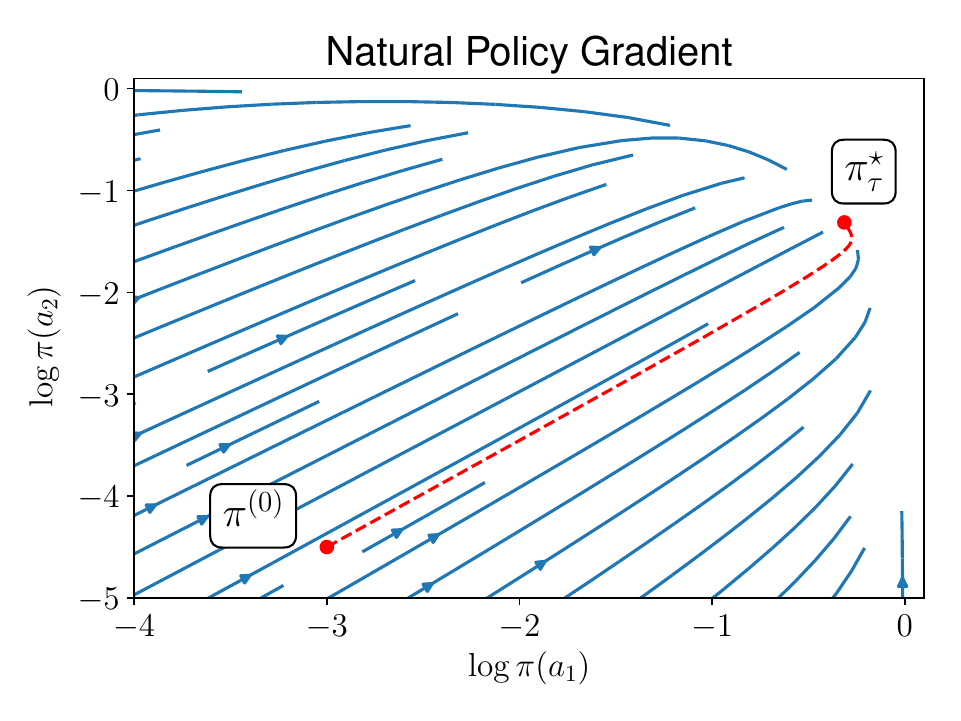}
    \end{tabular}
    \end{center}
    \caption{Comparison of PG and NPG methods with entropy regularization for a bandit problem $(\gamma = 0)$ with 3 actions associated with rewards $1.0$, $0.9$ and $0.1$. The regularization parameter is set to $\tau = 0.1$.}
\end{figure}

\section{Global convergence of policy gradient methods in games}
\label{sec:games}

In this section, we review recent progresses in developing policy gradient methods for games and multi-agent RL, focusing on the basic models of two-player zero-sum matrix games and two-player zero-sum Markov games.

\subsection{Problem settings}

\paragraph{Two-player zero-sum matrix games.} Given two players each taking actions from action spaces $\cA$ and $\cB$, the matrix game aims to solve the following saddle-point optimization problem
\begin{align}
    \max_{\mu \in \Delta(\cA)}\min_{\nu \in \Delta(\cB) } V^{\mu, \nu} := \mu^\top A \nu  ,
	\label{eq:max-min-problem-games}
\end{align}
where $A\in\mathbb{R}^{|\cA|\times |\cB|}$ denotes the payoff matrix with $\|A\|_\infty \le 1$, $\mu\in \Delta(\cA)$ and $\nu \in \Delta(\cB)$ stand for the mixed/randomized policies of each player, defined respectively as distributions over the probability simplex $\Delta(\cA)$ and $\Delta(\cB)$. Here, one player seeks to maximize the value function (i.e., the max player) while the other player wants to minimize it (i.e., the min player). 
 It is well-known since \citet{neumann1928theorie} that the max and min operators in \eqref{eq:max-min-problem-games} can be exchanged without affecting the solution. A pair of policies $(\mu^{\star},\nu^{\star})$ is said to be a {\em Nash equilibrium (NE)} of the matrix game if
\begin{align}
	V^{\mu^{\star}, \nu} \geq V^{\mu^{\star}, \nu^{\star }} \geq V^{\mu, \nu^{\star }} \qquad \forall (\mu, \nu) \in \Delta(\cA) \times \Delta(\cB).  
	\label{eq:defn-Nash-equiv-matrix}
\end{align}
In words, the NE corresponds to when both players play their best-response strategies against their opponents.

\paragraph{Two-player zero-sum Markov games.} Moving onto sequential decision-making, we consider an infinite-horizon discounted Markov game which is defined as $\mathcal{M} = \{\mathcal{S}, \cA, \cB, P, r, \gamma\}$, with discrete state space $\mathcal{S}$, action spaces of two players $\cA$ and $\cB$, transition probability $P$, reward function $r: \mathcal{S}\times\cA\times\cB \to [0, 1]$ and discount factor $\gamma \in [0, 1)$. 
A policy $\mu: \mathcal{S} \to \Delta (\cA)$ (resp. $\nu: \mathcal{S} \to \Delta (\cB)$) defines how the max player (resp. the min player) reacts to a given state $s$, where the probability of taking action $a \in \cA$ (resp. $b\in\cB$) is $\mu(a|s)$ (resp. $\nu(b|s)$). The transition probability kernel $P: \mathcal{S}\times \cA \times \cB \to \Delta(\mathcal{S})$ defines the dynamics of the Markov game, where $P(s'|s,a,b)$ specifies the probability of transiting to state $s'$ from state $s$ when the players take actions $a$ and $b$ respectively. The value function and $Q$-function of a given policy pair $(\mu, \nu)$ is defined in a similar way as in single-agent RL: 
\begin{align*}
    V^{\mu, \nu}(s) &= \mathbb{E} \left[ \sum_{t=0}^{\infty} \gamma^t r(s_t, a_t, b_t ) \,\big|\, s_0 =s \right],\\
    Q^{\mu, \nu}(s, a, b) &= \mathbb{E} \left[ \sum_{t=0}^{\infty} \gamma^t r(s_t, a_t, b_t ) \,\big|\, s_0 =s, a_0 = a, b_0 = b \right].
\end{align*}

The minimax game value on state $s$ is defined by
\[
    V^\star(s) = \max_{\mu}\min_{\nu}V^{\mu, \nu}(s) = \min_{\nu}\max_{\mu}V^{\mu, \nu}(s).
\]
Similarly, the minimax Q-function $Q^\star(s, a, b)$ is defined by
\begin{equation} 
\label{eq:Q_star_game}
Q^\star(s, a, b) = r(s, a, b) + \gamma\mathbb{E}_{s'\sim P(\cdot|s, a, b)}V^\star(s').
\end{equation}
It is established by \cite{shapley1953stochastic} that there exists a pair of stationary policy $(\mu^{\star}, \nu^{\star})$ attaining the minimax value on all states \citep{filar2012competitive}, and is called the NE of the Markov game. We seek to obtain a pair of $\epsilon$-optimal NE policy pair --- denoted by $\epsilon$-NE --- $(\hat{\mu}^\star, \hat{\nu}^\star)$ that satisfies
\begin{equation*}
    V^{\mu, \hat{\nu}^\star}(s) - \epsilon \le V^{\hat{\mu}^\star, \hat{\nu}^\star}(s) \le V^{\hat{\mu}^\star, \mu}(s) + \epsilon
\end{equation*}
for any $\mu \in \Delta(\cA)^\cS$, $\nu \in \Delta(\cB)^\cS$ and $s\in\mathcal{S}$.

With the success of the NPG method in single-agent RL, it is attempting to apply it in the context of two-player zero-sum games, where each player executes the NPG updates independently by treating the other player as part of the environment. 
Notably, the NPG dynamics coincide with the well-studied MWU method, or Hedge, that stems from online optimization and game theory, whose \textit{average-iterate} policy (i.e., $(\frac{1}{T}\sum_{t=1}^T\mu^{(t)}, \frac{1}{T}\sum_{t=1}^T\nu^{(t)})$) is shown to achieve a converging NE-gap at the rate of $\mathcal{O}(1/\sqrt{T})$ for two-player zero-sum matrix games. However, two technical challenges remain: last-iterate convergence and generalization to Markov games, which we shall discuss separately in the sequel.

\subsection{Global last-iterate convergence}

Average-iterate convergence guarantees fall short of determining if the learning trajectory converges towards NE (referred to as \textit{last-iterate} convergence that applies to $(\mu^{(T)}, \nu^{(T)})$ alone) or enters recurrent cycles instead. In addition, for large-scale applications that involve the use of neural networks, average-iterate convergence is also unsatisfactory as averaging neural networks is intractable. \citet{mertikopoulos2018cycles} demonstrated that MWU, when adopted by both agents in a two-player zero-sum matrix game, suffers from the Poincar\'e recurrence phenomenon that forbids the method from converging. This holds for other methods in the family of ``Follow the Regularized Leader'' (FTRL), which necessitates algorithmic modifications to the original NPG/MWU updates.

\begin{figure}
    \centering
    \includegraphics[trim={50 120 50 20},clip,width=8cm]{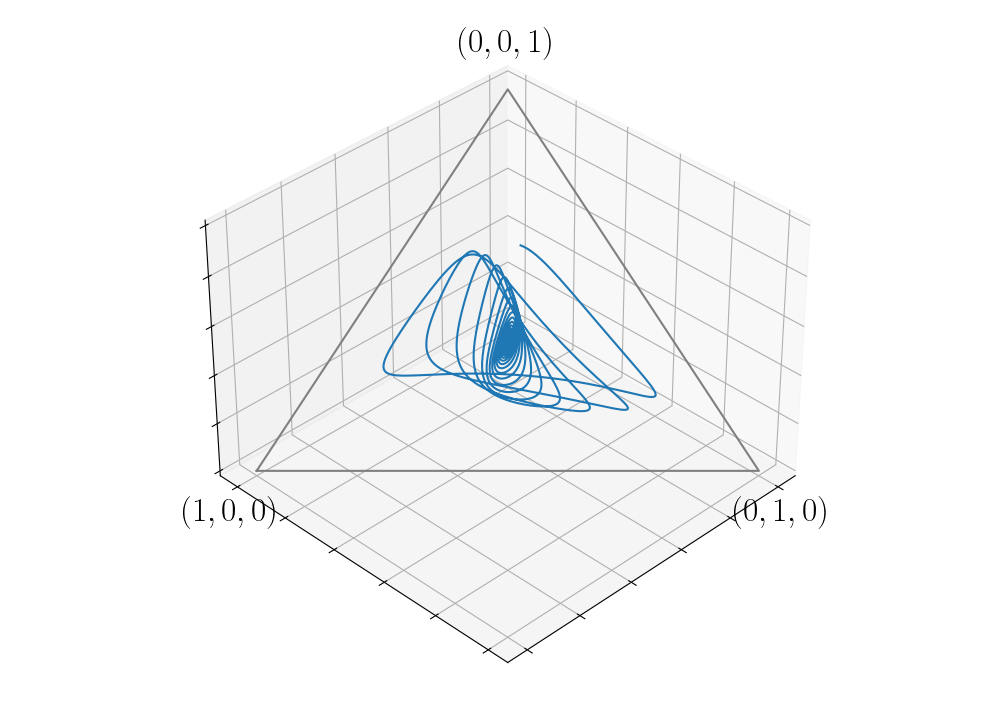}
    \caption{Cycles of MWU trajectories in a rock-paper-scissors game from a benign initialization close to the equilibrium.}
\end{figure}

\paragraph{Optimism.} \cite{rakhlin2013optimization} proposed the use of \textit{optimistic updates}, which extrapolate the gradient infomation with that from the previous iteration. In the context of two-player zero-sum matrix games, Optimistic MWU (OMWU) updates at the $t$-th iteration can be written as
\begin{equation*}
    \begin{cases}
        {\mu}^{(t+1)}(a) &\propto {\mu}^{(t)}(a) \exp(\eta (2A{\nu}^{(t)} - A{\nu}^{(t-1)}))\\
        {\nu}^{(t+1)}(b) &\propto {\nu}^{(t)}(b) \exp(-\eta (2A^\top{\mu}^{(t)} - A^\top{\mu}^{(t-1)}))
    \end{cases}.
\end{equation*}
Compared with the original MWU method
\begin{equation*}
    \begin{cases}
        {\mu}^{(t+1)}(a) &\propto {\mu}^{(t)}(a) \exp(\eta A{\nu}^{(t)})\\
        {\nu}^{(t+1)}(b) &\propto {\nu}^{(t)}(b) \exp(-\eta A^\top{\mu}^{(t)})
    \end{cases},
\end{equation*}
OMWU estimates the reward for the next iteration by appending the terms $A{\nu}^{(t)} - A{\nu}^{(t-1)}$ and $A^\top{\mu}^{(t)} - A^\top{\mu}^{(t-1)}$ to the current payoff vectors $A{\nu}^{(t)}$ and $A^\top{\mu}^{(t)}$. The updates can be equivalently formalized as Algorithm \ref{alg:OMWU}, where $\{\bar{\mu}^{(t)}, \bar{\nu}^{(t)}\}_{t=0}^\infty$ are auxiliary variables that can be viewed as some ``predictive'' sequence for facilitating the analysis.
\begin{algorithm}[h]
    \caption{OMWU for two-player zero-sum matrix games}
    \label{alg:OMWU}
    \textbf{Input:} Learning rate $\eta$, (optional) regularization parameter $\tau$.\\
    \textbf{Initialization:} Set $\mu^{(0)}$, $\nu^{(0)}$, $\bar\mu^{(0)}$ and $\bar\nu^{(0)}$ as uniform policies. Set $\tau=0$ when not using regularization.\\
    \For{$t=0,\cdots,\infty$}{
        When $t \ge 1$, update $\bar{\mu}^{(t)}$ and $\bar{\nu}^{(t)}$ as
        \begin{equation*}
            \begin{cases}
                \bar{\mu}^{(t)}(a) &\propto \bar{\mu}^{(t-1)}(a)^{1-\eta\tau} \exp(\eta A{\nu}^{(t)})\\
                \bar{\nu}^{(t)}(b) &\propto \bar{\nu}^{(t-1)}(b)^{1-\eta\tau} \exp(-\eta A^\top{\mu}^{(t)}).
            \end{cases}
            \end{equation*}
        Update ${\mu}^{(t+1)}$ and ${\nu}^{(t+1)}$ as
            \begin{equation*}
                \begin{cases}
                    {\mu}^{(t+1)}(a) &\propto \bar{\mu}^{(t)}(a)^{1-\eta\tau} \exp(\eta A{\nu}^{(t)})\\
                    {\nu}^{(t+1)}(b) &\propto \bar{\nu}^{(t)}(b)^{1-\eta\tau} \exp(-\eta A^\top{\mu}^{(t)})
                \end{cases}.
        \end{equation*}        
    }
\end{algorithm}
\citet{rakhlin2013optimization} first proved that OMWU yields $\cO(\log T)$ regret in two-player zero-sum matrix games, thus achieving a faster $\widetilde{\cO}(1/T)$ average-iterate convergence to NE. \citet{daskalakis2021near} demonstrated that OMWU achieves near-optimal $\cO(\text{poly}(\log T))$ regret for multi-player general-sum games as well. \citet{daskalakis2018last} established asymptotic last-iterate convergence of OMWU assuming that the NE is unique. \citet{wei2020linear} demonstrated OMWU converges linearly to the NE under the same uniqueness assumption.
\begin{theorem}[\mbox{\cite[informal]{wei2020linear}}]
    For a two-player zero-sum matrix game with unique NE $\zeta^\star = (\mu^\star, \nu^\star)$, the last iterate of OMWU $\zeta^{(T)} = (\mu^{(T)}, \nu^{(T)})$ converges to $\zeta^\star$ linearly with constant learning rate $\eta \le 1/8$.
\end{theorem}
\citet{wei2020linear} also investigated optimistic gradient descent ascent (OGDA), another variant of optimistic update rules by focusing on projected gradient updates \eqref{eq:direct_pg}, and derived global linear convergence guarantees without placing assumptions on NE uniqueness. A concrete iteration complexity, however, remains elusive as both convergence rates depend on unspecified problem-dependent parameters.

\paragraph{Regularization.} Regularization has been proven instrumental in enabling faster convergence for single-agent RL. In the context of game theory, entropy regularization is closely related to quantal response equilibrium (QRE) \citep{mckelvey1995quantal}, an extension to NE with bounded rationality. Formally speaking, when the two agents seek to maximize their own entropy-regularized payoffs 
\begin{equation*}
    \max_{\mu \in \Delta(\cA)}\min_{\nu \in \Delta(\cB) } V_\tau^{\mu, \nu} := \mu^\top A \nu + \tau \cH(\mu) - \tau \cH(\nu),
\end{equation*} 
the resulting equilibrium $\zeta_\tau^\star = (\mu_\tau^\star, \nu_\tau^\star)$ is referred to as a QRE and satisfies
\begin{equation*}
    \begin{cases}
        \mu_\tau^\star (a) &\propto \exp([A\nu_\tau^\star/\tau]_a) \\
        \nu_\tau^\star (b) &\propto \exp([A^\top\mu_\tau^\star/\tau]_b) \\
    \end{cases}, \qquad \forall a\in\cA, b\in\cB.
\end{equation*}
An approximate $\epsilon$-QRE is defined in a similar way to $\epsilon$-NE, by replacing $V^{\mu, \nu}$ in \eqref{eq:defn-Nash-equiv-matrix} with its regularized counterpart $V_\tau^{\mu, \nu}$, and $\epsilon/2$-QRE is guaranteed to be an $\epsilon$-NE by setting $\tau = \widetilde{\cO}(\epsilon)$. \citet{cen2021fast} proposed entropy-regularized OMWU (summarized in Algorithm \ref{alg:OMWU}) to combine the ideas of regularization and optimism, and proved the method converges to the QRE at a linear rate without assumption on the uniqueness of NE.
\begin{theorem}[\mbox{\cite[informal]{cen2021fast}}]
    With constant learning rate $0 < \eta \le \min\{1/(2\tau+2), 1/4\}$, the last iterate $\zeta^{(T)} = (\mu^{(T)}, \nu^{(T)})$ generated by entropy-regularized OMWU converges to the QRE $\zeta_\tau^\star$ at a linear rate $1-\eta\tau$.
\end{theorem}
This gives an iteration complexity of $\cO\big((1+1/\tau)\log1/\epsilon\big)$ for finding an $\epsilon$-QRE in a last-iterate manner, or $\widetilde{\cO}\big(1/\epsilon\big)$ for finding an $\epsilon$-NE by setting $\tau = \widetilde{\cO}(\epsilon)$. One might wonder if using regularization alone can also ensure the last-iterate convergence, which amounts to studying the update rule
\begin{equation*}
    \begin{cases}
        {\mu}^{(t)}(a) &\propto {\mu}^{(t-1)}(a)^{1-\eta\tau} \exp(\eta A{\nu}^{(t)})\\
        {\nu}^{(t)}(b) &\propto {\nu}^{(t-1)}(b)^{1-\eta\tau} \exp(-\eta A^\top{\mu}^{(t)})
    \end{cases}.
\end{equation*}
The answer turns out to be yes, as investigated recently in \citet{sokota2022unified,pattathil2023symmetric}. The same contraction rate $1-\eta\tau$ can be obtained albeit with a more restrictive choice of learning rate $\eta = \cO(\tau)$, which  translates to an iteration complexity of $\widetilde{\cO}(1/\epsilon^2)$ for finding an $\epsilon$-NE, slower than the entropy-regularized OMWU by a factor of $\widetilde{\cO}\big(1/\epsilon\big)$.  

\subsection{Extension to Markov games} 
Given some prescribed initial state distribution $\rho$, one can focus on solving the following saddle-point optimization problem
$$  \max_{\mu \in \Delta(\cA)^{|\cS|}}\min_{\nu \in \Delta(\cB)^{|\cS|} }  V^{\mu, \nu}(\rho) $$ 
for two-player zero-sum Markov games. A key property  for two-player zero-sum matrix games is that the value function $V^{\mu, \nu}$ is bilinear in the policy space, which unfortunately no longer holds in the Markov setting. Two strategies have been successful in establishing provable policy gradient methods for two-player zero-sum Markov games, by either leveraging tools from saddle-point optimization theory for nonconvex-nonconcave functions or exploiting recursive structures of the value function. 

Conventional wisdom in nonconvex-nonconcave saddle-point optimization suggests adopting two-timescale rules, which enforces a much smaller learning rate on one of the players.
\citet{daskalakis2020independent} demonstrated that when the two players adopt projected gradient descent ascent (GDA) updates with two-timescale learning rates, the method achieves an average-iterate convergence to $\epsilon$-NE within $\widetilde{\mathcal{O}}(\text{poly}(\epsilon^{-1},|\cS|,|\cA|,|\cB|))$ iterations (omitting additional instance-dependent parameters). \citet{zeng2022regularized} showed that softmax policy GDA updates with entropy regularization yields a last-iterate convergence with an improved iteration complexity with regard to $\epsilon^{-1}$. 

However, two-timescale rules give asymmetric convergence guarantees, i.e., only the slow learner is guaranteed to find an approximate NE policy, while the fast learner approximates the best response to the slow learner throughout the learning process. A natural question arises: is it possible to design symmetric algorithms with an improved iteration complexity?

\paragraph{Smooth value updates.} 
Instead of sticking to policy updates using vanilla gradient information $\nabla_\mu V^{\mu, \nu}(\rho)$ and $\nabla_\nu V^{\mu, \nu}(\rho)$, another line of works seek to divide the updates into two parts, where the policy updates are provided by two-player zero-sum matrix game algorithms with $Q^{(t)}(s) = [Q^{(t)}(s,a,b)]_{a\in\cA,b\in\cB}$, where
\begin{subequations} \label{eq:critic}
\begin{equation}
Q^{(t)}(s,a,b) =   r(s,a,b) + \gamma \ex{s'\sim P(\cdot|s,a,b)}{V^{(t-1)}(s')} ,
\end{equation} 
presuming the role of payoff matrices for all $s\in\cS$, and the value updates are given by
\begin{equation}    
        V^{(t)}(s) = (1-\alpha_{t}) V^{(t-1)}(s) + \alpha_{t} f^{(t)}(s).
\end{equation}
\end{subequations}
Here, $\alpha_{t} > 0$ is the learning rate for the value function and $f^{(t)}(s)$ is a one-step look-ahead value estimator for state $s$, typically defined as
\begin{equation}
    f^{(t)}(s) = \mu^{(t)}(s)^\top Q^{(t)}(s)\nu^{(t)}(s),
    \label{eq:f_unreg}
\end{equation}
or 
\begin{align}
    f^{(t)}(s) &= \mu^{(t)}(s)^\top Q^{(t)}(s)\nu^{(t)}(s)  + \tau \cH(\mu^{(t)}(s)) - \tau \cH(\nu^{(t)}(s))
    \label{eq:f_reg}
\end{align}
when we incorporate entropy regularization. These methods are akin to the prevalent actor-critic type algorithms. The algorithm procedure is summarized in Algorithm~\ref{alg:AC_Markov}. Note that when setting $\alpha_{t} = 1$ and $f^{(t)}(s)$ to the one-step minimax game value as
\begin{equation*}
    f^{(t)}(s) = \max_{\mu(s)\in\Delta(\cA)}\min_{\nu(s)\in\Delta(\cB)} \mu(s)^\top Q^{(t)}(s)\nu(s),
\end{equation*}
we recover the classical value iteration for two-player zero-sum Markov games.

\begin{algorithm}
    \caption{Actor-critic for two-player zero-sum Markov games}
    \label{alg:AC_Markov}
    \textbf{Input:} Learning rate for $Q$-value function $\{\alpha_t\}_{t=0}^{\infty}$, learning rate for policies $\eta$, policy optimization method for two-player zero-sum matrix game $\texttt{matrix\_alg}$, (optional) regularization parameter $\tau$.\\
    \textbf{Initialization:} Set $Q^{(0)} = 0$ and $\mu^{(0)}$, $\nu^{(0)}$ as uniform policies.\\
    \SetKwProg{ForP}{for}{ do in parallel}{end}
    \For{$t=0,1,\cdots$}{
        \ForP{all $s\in\cS$}{
            Invoke $\texttt{matrix\_alg}$ with payoff matrix $Q^{(t)}(s)$ and learning rate $\eta$ to update $\mu^{(t+1)}(s)$, $\nu^{(t+1)}(s)$.\\
            Update $Q^{(t+1)}(s)$ and $V^{(t+1)}(s)$ according to \eqref{eq:critic} with learning rate $\alpha_{t+1}$.
        }
    }
\end{algorithm}

\citet{wei2021last} first demonstrated that Algorithm~\ref{alg:AC_Markov} with $\texttt{matrix\_alg}$ OGDA and decaying learning rate $\alpha_t = \frac{2/(1-\gamma)+1}{2/(1-\gamma)+t}$ yields both average-iterate and last-iterate convergences to NE. \citet{cen2022faster} achieved an improved convergence rate by adopting entropy regularization and the OMWU method (e.g.~Algorithm~\ref{alg:OMWU}).
\begin{theorem}[\mbox{\cite[Theorem 1]{cen2022faster}}]
    Algorithm~\ref{alg:AC_Markov} with $\texttt{matrix\_alg}$ entropy-regularized OMWU, entropy-regularized value updates \eqref{eq:f_reg} and constant learning rates $\eta = \cO((1-\gamma)^3/|\cS|)$, $\alpha_t = \eta\tau$ guarantees last-iterate convergence to the QRE at a linear rate $1-\eta\tau$.
\end{theorem}
The above theorem demonstrates an iteration complexity of $\widetilde{\cO}\Big(\frac{|\cS|}{(1-\gamma)^4 \tau}\log1/\epsilon\Big)$ for finding an $\epsilon$-QRE, or $\widetilde{\cO}\Big(\frac{|\cS|}{(1-\gamma)^5 \epsilon}\Big)$ for finding an $\epsilon$-NE. It remains an open problem to achieve an iteration complexity with better dependency on $|\cS|$ and $(1-\gamma)^{-1}$.

\section{Global convergence of policy gradient methods in control}
\label{sec:control}
In this section, we briefly review policy gradient methods for control, focusing on a standard control problem called linear quadratic regulators (LQRs), based primarily on the excellent work of \citet{fazel2018global}. We refer interested readers to \citet{hu2023toward} for a recent comprehensive survey on the developments of policy optimization for general control problems including but not limited to $\mathcal{H}_\infty$ control, risk-sensitive control, to name just a few.

\subsection{Problem settings}

\paragraph{Linear quadratic regulator (LQR).}

Consider a discrete-time linear dynamic system
\begin{equation*}
x^{(t+1)} = A x^{(t)} + B u^{(t)},
\end{equation*}
where $x^{(t)}\in \mathbb{R}^{d}$ and $u^{(t)}\in \mathbb{R}^{k}$ are the state and the input at time $t$, $A\in\mathbb{R}^{d\times d}$ and $B \in \mathbb{R}^{d\times k}$ specify system transition matrices. The linear quadratic regulator (LQR) problem in the infinite horizon is defined as
\begin{equation}
\min_{u^{(t)}} \; \ex{x^{(0)}\sim \mathcal{D}}{\sum_{t=0}^\infty \Big(x^{(t)\top} Q x^{(t)} + u^{(t)\top} R u^{(t)}\Big)},
\end{equation}
where $\mathcal{D}$ determines the distribution of the initial state $x^{(0)}$, $Q \in \mathbb{R}^{d\times d}$ and $R\in \mathbb{R}^{k\times k}$ are positive definite matrices parametrizing the costs. Classical optimal control theory \citep{anderson2007optimal} tells us that under certain stability conditions (e.g., controllability), it is ensured that the optimal cost is finite and can be achieved by a linear controller
\begin{equation*}
    u^{(t)} = - K^\star x^{(t)},
\end{equation*}
where $K^{\star}\in\mathbb{R}^{k\times d}$ is the optimal control gain matrix. From an optimization perspective,
it is therefore natural to cast the LQR problem as optimizing over all linear controllers $u^{(t)} = -Kx^{(t)}$ with $K\in \mathbb{R}^{k\times d}$ to minimize the cost:
\begin{equation*}
    \min_{K\in\mathcal{K}}\; C(K) := \ex{x^{(0)}\sim \mathcal{D}}{x^{(0)\top} P_K x^{(0)}},
\end{equation*}
where 
\begin{equation*}
    P_K = \sum_{t=0}^{\infty} ((A - BK)^\top)^t (Q+K^\top RK)(A-BK)^t.
\end{equation*}
The feasible set $\mathcal{K} = \{K: \|A-BK\|_2< 1\}$ ensures $P_K$ is well defined for all $K \in \mathcal{K}$.

\subsection{Policy gradient method}
The policy gradient method for the LQR problem is simply defined as
\begin{equation}
    K^{(t+1)} = K^{(t)} - \eta \nabla_K C(K^{(t)}),
    \label{eq:pg_lqr}
\end{equation}
where $\eta>0$ is the learning rate. In addition, the policy gradient $\nabla_K C(K)$ can be written as
\begin{equation*}
    \nabla_K C(K) = 2\big((R + B^\top P_K B)K - B^\top P_K A \big) \Sigma_K,
\end{equation*}
where $\Sigma_K=  \ex{x^{(0)}\sim \mathcal{D}}{\sum_{t=0}^\infty x^{(t)}x^{(t)\top}}$ denotes the state correlation matrix.

Like the RL problem, the objective function $C(K)$ is nonconvex in general, which makes it challenging to claim a global convergence guarantee. Fortunately, the LQR problem satisfies the following gradient dominance condition \citep{fazel2018global}.
\begin{lemma}[Gradient dominance]
    Suppose that $ \lambda = \sigma_{\min}(\ex{x^{(0)}\sim \mathcal{D}}{x^{(0)}x^{(0)\top}}) > 0$. It holds that
    \begin{equation*}
        C(K) - C(K^\star) \le \frac{\norm{\Sigma_{K^\star}}_2}{\lambda^2 \sigma_{\min}(R)}\norm{\nabla_K C(K)}_{\mathsf{F}}^2.
    \end{equation*}
\end{lemma}
The gradient dominance property provides hope for attaining global linear convergence of the policy gradient method to the optimal policy, yet to complete the puzzle, another desirable ingredient is the smoothness of the cost function $C(K)$. While this cannot be established in its full generality as $C(K)$ becomes infinity when $K$ moves beyond $\mathcal{K}$, fortunately, smoothness can be established for any given sublevel set $\mathcal{K}_{\bar{\gamma}} = \{K \in \mathcal{K}: C(K) \le \bar{\gamma}\}$, which suffices 
 to establish the desired convergence results as follows.
\begin{theorem}[\mbox{\citet{fazel2018global}}]
    Assume that $C(K^{(0)})$ is finite.
    With an appropriate constant learning rate 
    \begin{equation*}
        \eta = \text{poly}\Big(\frac{\lambda \sigma_{\min}(Q)}{C(K^{(0)})}, \norm{A}_2^{-1}, \norm{B}_2^{-1}, \norm{R}_2^{-1}, \sigma_{\min}(R)\Big),
    \end{equation*}
    the policy gradient method \eqref{eq:pg_lqr} converges to the optimal policy at a linear rate:
    \begin{equation*}
        C(K^{(t+1)}) - C(K^\star) \le \Big(1 - \frac{\lambda^2\sigma_{\min}(R)\eta}{\norm{\Sigma_{K^\star}}_2}\Big) \big(C(K^{(t)}) - C(K^\star) \big).
    \end{equation*}
\end{theorem}
To find  an $\epsilon$-optimal control policy $K^{(T)}$ satisfying $C(K^{(T)}) - C(K^\star) \le \epsilon$, the above theorem ensures that     the policy gradient method takes no more than
\begin{equation*}
     \frac{\norm{\Sigma_{K^\star}}_2}{\lambda^2 \sigma_{\min}(R)\eta}\log\frac{C(K^{(0)}) - C(K^\star)}{\epsilon}
\end{equation*}
iterations.

\subsection{Natural policy gradient method}
To facilitate the development of NPG, we consider a linear policy with additive Gaussian noise, specified as
\begin{equation*}
    u^{(t)} \sim \pi(\cdot|x^{(t)}) = \mathcal{N}( -Kx^{(t)}, \sigma^2I).
\end{equation*}
The NPG update rule \citep{kakade2002natural} then reads like
\begin{equation}
    \text{vec}(K^{(t+1)}) = \text{vec}(K^{(t)}) - \eta (\mathcal{F}_{\mathcal{D}}^{K^{(t)}})^\dagger \text{vec}(\nabla_K C(K^{(t)})),
    \label{eq:control_NPG}
\end{equation}
where $\text{vec}(K)$ flattens $K \in \mathbb{R}^{k\times d}$ into a vector in row-major order, and the Fisher information matrix $\mathcal{F}_{\mathcal{D}}^{K}\in\mathbb{R}^{kd\times kd}$  is given by
\begin{align*}
    \mathcal{F}_{\mathcal{D}}^K &= \ex{}{\sum_{t=0}^\infty \text{vec}(\nabla_K \log\pi(u^{(t)}|x^{(t)}))\text{vec}(\nabla_K \log\pi(u^{(t)}|x^{(t)}))^\top}\\
    &= \sigma^{-2} \ex{}{\sum_{t=0}^\infty \text{diag}(x^{(t)}x^{(t)\top},\cdots, x^{(t)}x^{(t)\top})}\\
    &= \sigma^{-2} \text{diag}(\Sigma_{K},\cdots,\Sigma_{K}).
\end{align*}
Merging the dummy variance $\sigma^2$ into the learning rate $\eta$, and reshaping back into the matrix form, the NPG update rule \eqref{eq:control_NPG} can be equivalently rewritten as
\begin{equation*}
    K^{(t+1)} = K^{(t)} - \eta \nabla_K C(K^{(t)})\Sigma_{K^{(t)}}^{-1},
\end{equation*}
which modifies the update direction using the state correlation matrix. This allows for an improved progress following a single update, as demonstrated by the following lemma.
\begin{lemma}[\mbox{\citet{fazel2018global}}]
    Assume that $C(K^{(t)})$ is finite and the learning rate satisfies $\eta \le \frac{1}{\|R + B^\top P_{K^{(t)}} B\|_2}$, the NPG update satisfies
    \begin{equation*}
        C(K^{(t+1)}) - C(K^\star) \le \Big(1- \frac{\lambda  \sigma_{\min}(R) \eta}{\|\Sigma_{K^\star}\|_2}\Big)\big(C(K^{(t)}) - C(K^\star)\big).
    \end{equation*}
\end{lemma}
The improvement is twofold: the convergence rate is improved by a factor of $\lambda$; in addition, it allows for a larger learning rate that would not be possible under a smoothness-based analysis. By adopting $\eta = \frac{1}{\|R\|_2 + \|B\|_2^2C(K^{(0)})\lambda^{-1}}$, one can show that the learning rate requirement is met throughout the trajectory, which leads to an iteration complexity of 
\begin{equation*}
  \frac{\norm{\Sigma_{K^\star}}_2}{\lambda \sigma_{\min}(R)\eta}\log\frac{C(K^{(0)}) - C(K^\star)}{\epsilon}
\end{equation*}
for finding an $\epsilon$-optimal control policy. Last but not least, it is possible to achieve even-faster convergence rate by assuming access to more complex oracles, i.e., the Gauss-Newton method. The update rule is given by
\begin{equation*}
    K^{(t+1)} = K^{(t)} - \eta (R + B^\top P_{K^{(t)}}B)^{-1} \nabla_K C(K^{(t)})\Sigma_{K^{(t)}}^{-1},
\end{equation*}
which allows a constant learning rate as large as $\eta = 1$ and an improved iteration complexity of $$  \frac{\norm{\Sigma_{K^\star}}_2}{\lambda}\log\frac{C(K^{(0)}) - C(K^\star)}{\epsilon}. $$

\section{Conclusions}
\label{sec:conclusions}

Policy gradient methods remain to be at the forefront of data-driven sequential decision making, due to its simplicity and flexibility in integrating with other advances in computation  from adjacent fields such as high performance computing and deep learning. Due to space limits, our focus is constrained on the optimization aspects of policy gradient methods, assuming access to exact gradient or policy evaluations. In reality, these information need to estimated by samples collected via various mechanisms and thus noisy, leading to deep interplays between statistics and optimization. We hope this exposition provides a teaser to invite more interest from the optimization community to work in  this area.

\paragraph*{Acknowledgement}
S.~Cen and Y. Chi are supported in part by the grants ONR N00014-19-1-2404, NSF CCF-1901199, CCF-2106778 and CNS-2148212. S.~Cen is also gratefully supported by Wei Shen and Xuehong Zhang Presidential Fellowship, Boeing Scholarship, and JP Morgan Chase AI PhD Fellowship.

\bibliographystyle{apalike} 

\bibliography{bibfileRL.bib,bibfileGame.bib,bibfileControl.bib}

\end{document}